\documentclass[11pt]{amsart}
\usepackage{amssymb}
\usepackage{amsfonts}
\usepackage{amscd}
\usepackage{amsmath,amsfonts,amsthm,amscd,amssymb,euscript,
graphics,psfrag}

\newcommand{\beq}{\begin{equation}}
\newcommand{\eeq}{\end{equation}}

\newtheorem{theorem}{Theorem}[section]
\newtheorem{lemma}{Lemma}[section]
\newtheorem{corollary}{Corollary}[section]

\theoremstyle{definition}

\newtheorem*{ip}{Independence Problem}

\input tcilatex

\begin{document}
\date{October 2009}
\title{Strong uniform expansion in $\mathrm{SL}(2,p)$}
\author{Emmanuel Breuillard}
\address{Laboratoire de Math\'{e}matiques, Universit\'{e} Paris-Sud 11,
91405 Orsay cedex, France}
\email{Emmanuel.Breuillard@math.u-psud.fr}
\author{Alex Gamburd}
\address{Department of Mathematics, University of California at Santa Cruz,
1156 High Street, Santa Cruz, CA 95064, USA}
\email{agamburd@ucsc.edu}

\begin{abstract}
We show that there is an infinite set of primes $\mathcal{P}$ of
density
one, such that the family of \textit{all} Cayley graphs of $\mathrm{SL}(2,p)$%
, $p\in \mathcal{P}$, is a family of expanders.
\end{abstract}

\thanks{The first author was supported in part by ERC starting grant. The
second author was supported in part by DARPA, NSF, and the Sloan
Foundation.} \maketitle

\section{Introduction and statement of results}

Expanders are highly-connected sparse graphs widely used in
Computer Science; they also have found some remarkable
applications in pure
mathematics \cite{hlw, Sa04}. Given an undirected $d$-regular graph $%
\mathcal{G}$ and a subset $X$ of $V$, the \emph{expansion} of $X$,
$c(X)$, is defined to be the ratio $|\partial(X)|/|X|$, where
$\partial(X) = \{y \in \mathcal{G} \, : \, \mathrm{distance}(y, X)
= 1\} $. The \emph{expansion coefficient} of a graph $\mathcal{G}$
is defined as follows:
\begin{equation*}
c(\mathcal{G})=\inf\left \{c(X) \mid |X| < \frac{1}{2}
|\mathcal{G}|\right\}.
\end{equation*}

A family of graphs $\mathcal{G}_{n}$ forms a family of
$c$-expanders if there is $c>0$, such that
\begin{equation}
\inf_{n\in \mathbb{N}}c(\mathcal{G}_{n})\geq c. \label{e:exp}
\end{equation}

Usually one takes the family of graphs to be $d$-regular, in which
case the condition \eqref{e:exp} has an alternative spectral
interpretation. The \emph{adjacency matrix} of $\mathcal{G}$,
$A(\mathcal{G})$ is the $\mid \! \mathcal{G} \!\mid$ by $\mid \!
\mathcal{G} \!\mid$ matrix, with rows and columns indexed by
vertices of $\mathcal{G}$, such that the $x,y$ entry is 1 if and
only if $x$ and $y$ are adjacent and 0 otherwise.
Using the discrete analogue of Cheeger-Buser inequality, proved by
Alon and Milman, condition (\ref{e:exp}) can be rewritten in terms
of the second
largest eigenvalue 
of the adjacency matrix $A(\mathcal{G})$ as follows:
\begin{equation} \label{e:a1}
\sup_{n\in \mathbb{N}} \lambda_{1}(A_{n,d}) < d.
\end{equation}

Given a finite group $G$ with a symmetric set of generators $S$,
the Cayley graph $\mathcal{G}(G, S)$, is a graph which has
elements of $G$ as vertices and which has an edge from $x$ to $y$
if and only if $x = \sigma y$ for some $\sigma \in S$. The
explicit constructions of expander graphs
(by Margulis~\cite{Ma73, Ma88} and Lubotzky, Phillips and Sarnak~\cite{LPS86}%
) used deep tools (Kazhdan's property (T), Selberg's Theorem,
proved Ramanujan conjectures) to construct families of Cayley
graphs of finite
groups as follows. Starting with an infinite group $\Gamma$ (e.g. $\mathrm{SL%
}_2({\mathbb{Z}})$) and a finite set of generators $S$, $\langle S
\rangle =
\Gamma$, one considers a family of Cayley graphs $\mathcal{G}_{i}=\mathcal{G}%
(G_i, S_i)$, where $G_i $ is an infinite family of finite quotients (e.g. $%
\mathrm{SL}_2({\mathbb{F}_p}))$) and $S_i$ is an image of $S$
under the natural projection. On the other hand, as shown in
\cite{LW93}, some families of groups, for example abelian groups
or solvable groups of bounded derived length, cannot be made into
families of expanders with respect to \emph{any} choice of
generators. A basic problem formulated by Lubotzky and Weiss in
1993 \cite{LW93}, (see also \cite{AL94, AL95}) is to what extent
the expansion property is the property of the family of groups
$\{G_i\}$ alone, independent of the choice of generators:

\begin{ip}[\protect\cite{LW93}]
\label{ip} Let $\{G_i\}$ be a family of finite groups, $\langle
S_i \rangle = \langle S_{i}^{\prime }\rangle = G_i$ and $|S_i| <k$
, $|S_{i}^{\prime }|<k $. Does the fact that $\{\mathcal{G}(G_i,
S_i)\}$ is an expander family imply the same for
$\{\mathcal{G}(G_i, S_{i}^{\prime })\}$?
\end{ip}

It turned out that in general the answer is negative. In 2001
Alon, Lubotzky and Wigderson \cite{ALW01}, using the notion of
zig-zag product introduced in the paper of Reingold, Vadhan and
Wigderson \cite{RVW02}, constructed a family of groups $G_{i}$
which are expanders with respect to one choice of generators and
not with respect to another such choice. The groups $G_{i}$
are of the form $A_{i}\rtimes B_{i}$ where $B_{i}=\mathrm{SL}_{2}({\mathbb{F}%
_{p_{i}}})$ and $A_{i}=\mathbb{F}_{2}^{P_{i}}$, with $P_{i}={\mathbb{F}%
_{p_{i}}}\cup \{\infty \}$ for an infinite family of primes. In
another breakthrough, Kassabov \cite{kas1} proved that symmetric
groups can be made expanders with respect to explicit sets of
generators --- it is easy to see that symmetric groups are not
expanding with respect to $\{(12),(1,2,\ldots ,n)\}$.

In \cite{BG} it was proved that Cayley graphs of $\mathrm{SL}_2({\mathbb{F}_p%
})$ are expanders with respect to the projection of any fixed elements in $%
\mathrm{SL}(2, {\mathbb{Z}})$ generating a non-elementary
subgroup, and with respect to generators chosen at random in
$\mathrm{SL}_2({\mathbb{F}_p})$. In this note we prove that Cayley
graphs of $\mathrm{SL}_2({\mathbb{F}_p})$ are expanders with
respect to all generators for infinitely many primes, thus
obtaining a first example of an affirmative answer to
Lubotzky-Weiss problem.

For a finite group $G$ generated by $2$ elements let $g(G)$ be the
smallest spectral gap of the averaging operator
$\frac{1}{4}(a+a^{-1}+b+b^{-1})$ among all possible choices of a
pair $(a,b)$ in $G$ which generates the
group. We say that $G$ has uniform spectral gap at least $g>0$ if \ $%
g(G)\geq g$. We say that an infinite sequence $\{G_{n}\}_{n}$ of finite $2$%
-generated groups has \textit{strong uniform expansion} if $%
\inf_{n}g(G_{n})>0.$ For $k$-generated Cayley graphs with fixed
$k$ a uniform lower bound on the spectral gap is equivalent to
\eqref{e:a1}.

\begin{theorem}
\label{main}There is a function $\varepsilon (\delta )>0$ with
$\varepsilon (\delta )\rightarrow 0$ as $\delta \rightarrow 0$
such that for all $\delta
>0$ and all $X>1$ the number of rational primes $p$ less than $X$ for which $%
\mathrm{SL}(2,p)$ has uniform spectral gap less than $\delta $ is at most $%
X^{\varepsilon (\delta )}.$
\end{theorem}

We show furthermore that in our case the expansion property
remains uniform as the number of generators increases, thus
obtaining a result valid for all Cayley graphs regardless of the
number of generators :

\begin{corollary}
\label{c0} There is a function $\varepsilon (c)>0$ with
$\varepsilon (c)\rightarrow 0$ as $c\rightarrow 0$ such that for
all $c>0$ and all $X>1$ the number of rational primes $p$ less
than $X$ for which some Cayley graph of $\mathrm{SL}(2,p)$ fails
to be a $c$-expander is at most $X^{\varepsilon (c)}.$
\end{corollary}

According to the prime number theorem there are roughly $X/\log
(X)$ primes less that $X,$ hence Theorem \ref{main} produces an
infinite sequence of finite groups with strong uniform expansion.
A well-known result about the distribution of primes (Hoheisel's
theorem \cite{Ho}, see \cite{KI} \S 10.5) says that there is a
constant $\beta _{0}>0 $ such that for all large $X,$ there is at
least one prime between $X\ $and $X+X^{\beta _{0}}$ (one can take
any $\beta _{0}>\frac{7}{12}$ by \cite{Hu}). Hence the following
immediate consequence :

\begin{corollary}
\label{c1} For any $\beta \in (\beta _{0},1)$, there is a constant $%
c=c(\beta )>0$ and an infinite sequence of primes $p_{n}$ with
$p_{n+1}\leq
p_{n}+p_{n}^{\beta }$ such that for every $n$ every Cayley graph of $\mathrm{%
SL}(2,p_{n})$ is a $c$-expander.
\end{corollary}

We note that while a lower bound for the constant $c$ in the last
two corollaries can be computed effectively, our method does not
produce any explicit infinite family of primes for which Corollary
\ref{c1} is true.

Combining Corollary \ref{c1} with the main result of \cite{GP}, we
obtain an application to product replacement graphs \cite{pak}.
Given a group $G$, the product replacement graph $\Gamma _{k}(G)$
introduced in \cite{c+} in connection with computing in finite
groups is defined as follows. The
vertices of $\Gamma _{k}(G)$ consist of all $k$-tuples of generators $%
(g_{1},\dots ,g_{k})$ of the group $G$. For every $(i,j)$, $1\leq
i,j\leq
k,i\neq j$ there is an edge corresponding to transformations \thinspace\ $%
L_{i,j}^{\pm }$ \thinspace\ and \thinspace\ $R_{i,j}^{\pm }$ :
\begin{align*}
& R_{i,j}^{\pm }\,:\ (g_{1},\dots ,g_{i},\dots ,g_{k})\rightarrow
(g_{1},\dots ,g_{i}\cdot g_{j}^{\pm 1},\dots ,g_{k}), \\
& L_{i,j}^{\pm }\,:\ (g_{1},\dots ,g_{i},\dots ,g_{k})\rightarrow
(g_{1},\dots ,g_{j}^{\pm 1}\cdot g_{i},\dots ,g_{k}).
\end{align*}

The graphs \thinspace\ $\Gamma _{k}(G)$ \thinspace\ are regular,
of degree \thinspace\ $4\,k\,(k-1)$, possibly with loops and
multiple edges. Connectivity of $\Gamma _{k}(G)$ has been the
subject of intensive recent investigations; for
$G=\mathrm{SL}_{2}(p)$ and $k\geq 3$ it was established by Gilman
in \cite{gil}.

In the case of the free group $F_{k}$ the moves $L_{i,j}^{\pm }$ and $%
R_{i,j}^{\pm }$ defined above correspond to Nielsen moves on
$\Gamma _{k}(F_{k})$. For every group $G$, the set $\Gamma
_{k}(G)$ can be
identified with $E=\mathrm{Epi}(F_{k},G)$, the set of epimorphisms from $%
F_{k}$ onto $G$, and the group $A=\mathrm{Aut}(F_{k})$ acts on $E$
in the following way: if $\alpha \in A$ and $\varphi \in E$,
$\alpha (\varphi
)=\varphi \cdot \alpha ^{-1}$. A long-standing problem is whether $\mathrm{%
Aut}(F_{k})$ has property T for $k\geq 4$; in \cite{LP} Lubotzky
and Pak
observed that a positive answer to this problem implies the expansion of $%
\Gamma _{k}(G)$ for all $G$ and proved that $\Gamma _{k}(G)$ are
expanders when $G$ is nilpotent of class $l$ and both $k$ and $l$
are fixed.

In \cite{GP} the second author and Pak established a connection
between the expansion coefficient of the product replacement graph
$\Gamma_k(G)$ and the minimal expansion coefficient of a Cayley
graph of $G$ with $k$ generators,
and in particular showed that the product replacement graphs \, $\Gamma_k(%
\mathrm{SL}(2,p))$ \, form an expander family under assumption of
strong
uniform expansion of $\mathrm{SL}(2,p)$ on $k$ generators. Corollary \ref%
{cpr} is an immediate consequence of Corollary \ref{c1} and
Corollary 2 in \cite{GP}.

\begin{corollary}
\label{cpr} Let $k\geq 4.$ The family of product replacement graphs $%
\{\Gamma _{k}(\mathrm{SL}(2,p_{n}))\}_{n}$ forms a family of
expanders.
\end{corollary}

\section{Proofs}

First a few words about the strategy of proof. Roughly speaking Theorem \ref%
{main} follows from the combination of the results of the second
named
author with Jean Bourgain on the spectral gap for $\mathrm{SL}(2,p)$ \cite%
{BG} with the Strong Tits Alternative proved by the first named author \cite%
{B} together with a combinatorial argument based on the effective
arithmetic Nullstellensatz which we explain here. We now give
details.

\subsection{Strong Tits.}

Let $g(p)=g(\mathrm{SL}(2,p)).$ We will in fact prove the
following equivalent version of Theorem \ref{main}.

\begin{theorem}[reformulation of the main theorem]
\label{mainbis}For every $\varepsilon >0$ and every $A>1$ there is
$\delta
>0 $ such that for every $X>1$, the number of primes $p$ in the interval $%
[X,X^{A}]$ with $g(p)<\delta $ is at most $X^{\varepsilon }.$
\end{theorem}

Recall the statement of the Strong Tits Alternative proved in
\cite{B} (see also \cite{B2} for a proof in the special case of
$GL(2)$, which is enough for the purpose of this paper).

\begin{theorem}[\protect\cite{B}]
\label{ST}There is a universal constant $N$ such that any finite
symmetric set $S$ in $GL(2,\overline{\mathbb{Q}})$ which does not
generate a virtually solvable group has the property that some
words $w_{1}$ and $w_{2}$ of length at most $N$ in the elements of
$S$ will generate a free subgroup.
\end{theorem}

Consider the set $\mathcal{C}_{n}$ of all assignments which assign
$4$ paths of length $n$ starting at the identity in the free group
$F_{2}$ to every
pair $(w_{1},w_{2})\in B(N)^{2}$ ($B(N)$ is the ball of radius $N$ in $F_{2}$%
). Observe that $|\mathcal{C}_{n}|=4^{4Kn}$ where $K:=|B(N)^{2}|.$
Among those, consider the subset $\mathcal{D}_{n}\subset
\mathcal{C}_{n}$ made of
assignments all of whose $4$ paths $W_{1},...,W_{4}$ satisfy $%
[[W_{1},W_{2}],[W_{3},W_{4}]]\neq 1$ (i.e. such that the
associated reduced word is non trivial).

\begin{lemma}
Let $S_{n}^{(i)}$ for $i=1,...,4$ be $4$ independent simple random
walks on
the free group $F_{2}.$ Then $\mathbb{P}%
([[S_{n}^{(1)},S_{n}^{(2)}],[S_{n}^{(3)},S_{n}^{(4)}]]=1)\leq
e^{-\kappa n}$ for some explicit $\kappa >0.$
\end{lemma}

\textit{Proof.} By Kesten's theorem \cite{HK59} for every $x\in
F_{2},$ the probability that the simple random walk on $F_{2}$
starting at $1 $ visits $x $ at time $n$ is at most $\left(
\frac{\sqrt{3}}{2}\right) ^{n}. $ Observe also that the
centralizer of a non trivial element in $F_{2}$ is a cyclic
subgroup of $F_{2}$ and that any cyclic subgroup intersects
$B(n)\backslash
\{1\}$ in at most $2n$ elements. Thus if $[S_{n}^{(1)},x]=1$ for some given $%
x\neq 1,$ then $S_{n}^{(1)}$ may take only $2n+1$ possible values
in $B(n).$ We may now write
\begin{equation*}
\mathbb{P}([S_{n}^{(1)},S_{n}^{(2)}]=1)\leq \max_{x\in B(n)\backslash \{1\}}%
\mathbb{P}([S_{n}^{(1)},x]=1)+\mathbb{P}(S_{n}^{(2)}=1)\leq
(2n+2)\cdot \left( \frac{\sqrt{3}}{2}\right) ^{n}
\end{equation*}%
Furthermore, note that if $a$ and $b$ are fixed and not $1$, then
the set of $x$'s such that $[a,x]=b$ coincides, if non empty, with
a coset of the
centralizer of $a.$ In particular this set can intersect $B(n)$ in at most $%
4n+1$ elements. Using this we can now write:%
\begin{eqnarray*}
\mathbb{P}([[S_{n}^{(1)},S_{n}^{(2)}],[S_{n}^{(3)},S_{n}^{(4)}]]
&=&1)\leq
\max_{u\neq 1,a\in B(n)\backslash \{1\}}\mathbb{P}([u,[a,S_{n}^{(4)}]]=1)+%
\mathbb{P}([S_{n}^{(1)},S_{n}^{(2)}]=1)+\mathbb{P}(S_{n}^{(3)}=1) \\
&\leq &(8n+1)\cdot \max_{a,b\neq
1}\mathbb{P}([a,S_{n}^{(4)}]=b)+(2n+3)\cdot
\left( \frac{\sqrt{3}}{2}\right) ^{n} \\
&\leq &[(8n+1)(4n+1)+(2n+3)]\cdot \left( \frac{\sqrt{3}}{2}\right)
^{n}\leq e^{-\kappa n}
\end{eqnarray*}%
for some explicit $\kappa >0.$

\begin{corollary}
\label{c vs d}There is an explicit number $\alpha >0$ such that $|\mathcal{C}%
_{n}\backslash \mathcal{D}_{n}|\leq |\mathcal{C}_{n}|^{1-\alpha
}.$
\end{corollary}

Theorem \ref{ST} now implies that for every assignment $c\in
\mathcal{D}_{n}$ the algebraic subvariety $\mathcal{W}_{c}$ of
$(GL_{2}(\mathbb{C}))^{2}$ defined by the vanishing of the
corresponding $4$-fold commutators is contained in the subvariety
$\mathcal{V}_{sol}$ in $(GL_{2}(\mathbb{C}))^{2}$ of pairs which
generate a virtually solvable subgroup, because each
commutator would give a non trivial relation for each of the pairs $%
(w_{1},w_{2})\in B(N)^{2}.$ The subvariety $\mathcal{V}_{sol}$
coincides with the set of pairs that leave invariant a finite
subset of at most $M$ points on the projective line (for some
constant $M$). Both varieties are defined over $\mathbb{Z}$ (see
\cite{B2} \S 9 for more on this translation).

\subsection{Effective Nullstellensatz}

Observe further that $\mathcal{W}_{c}$ is defined by the vanishing
of $K$ words of length at most $4Nn.$ Here $K$ and $N$ are fixed
constants and $n$ will grow. If $P_{1}(c),...,P_{4K}(c)$ denote
the polynomials (with integer coefficients) in the $8$ standard
variables of $(GL_{2})^{2}$ (i.e. the
matrix entries) whose vanishing define $\mathcal{W}_{c}$, and $%
Q_{1},...,Q_{s}$ the polynomials (also with integer coefficients) defining $%
\mathcal{V}_{sol},$ then the Nullstellensatz asserts the existence
of polynomials $g_{ij}$ with integer coefficients and natural
numbers $a_{i}(c)$ and $e_{i}(c)$ such that for each $i=1,...,s$
\begin{equation*}
a_{i}(c)Q_{i}^{e_{i}(c)}=\sum_{j}g_{ij}P_{j}(c)
\end{equation*}%
Moreover, standard versions of the effective arithmetic
Nullstellensatz, based for instance on the classical Hermann
method such as in Masser and Wusholtz's paper (\cite{MW}, chapter
4), give bounds on the $a_{i}$'s (and on the other parameters too,
but we will only need the bound for $a_{i}$). In our context, the
polynomials $P_{j}(c)$ have degree at most $4Nn$ and height (i.e.
maximum modulus of coefficients) at most $H^{n}$ for some constant
$H$ (see e.g. \cite{B2} \S 9). Then the Nullstellensatz bounds
from \cite{MW} (Theorem 4.1.IV) give the existence of constants
$C,r\geq 1$ such that $a_{i}(c)\leq e^{Cn^{r}}.$

We will denote by $a(c)$ the product $a_{1}(c)\cdot ...\cdot
a_{s}(c),$ which is again bounded above by $e^{Cn^{r}}$ (for some
other constant $C,$ since $s$ is a constant).
\begin{equation}
a(c)\leq e^{Cn^{r}} \label{bound}
\end{equation}

We need the following lemma.

\begin{lemma}
\label{gen}There is $p_{0}>0$ such that for all primes $p>p_{0}$
we have : for all $(a,b)\in \mathrm{SL}(2,p)^{2}$ if
$Q_{i}(a,b)=0$ for all $i=1,...,s$ then $(a,b)$ does not generate
$\mathrm{SL}(2,p).$
\end{lemma}

\textit{Proof.} \ This follows from \cite{B2} \S 9 : the variety $\mathcal{V}%
_{sol}$ coincides with the $\mathbb{Z}$-scheme of pairs $(a,b)$
that leave
invariant a finite set of at most $M$ points on the projective line $\mathbb{%
P}^{1}$ (for some constant $M$). This condition is given by the
vanishing of one of finitely many resultant polynomials, the
product of which must lie in
the ideal generated by the $Q_{i}$'s. Thus if the $Q_{i}(a,b)$ vanish in $%
\mathbb{F}_{p}$, for $p$ large enough, the group generated by
$(a,b)\in \mathrm{SL}(2,p)^{2}$ must have a subgroup of index at
most $M$ which fixes a point in
$\mathbb{P}^{1}(\overline{\mathbb{F}_{p}})$, hence does not
generate $\mathrm{SL}(2,p)$ when $p$ is $>M$ say. \hfill $\square
$

Thus if $p$ is a prime not dividing $a(c),$ then for every generating pair $%
(a,b)\in \mathrm{SL}(2,p)$ there must be some $j$ such that
$P_{j}(c)\neq 0,$ i.e. there must be 4 paths of length $n$ and two
words $w_{1},w_{2}$ in $a,b$ such that the resulting commutator
word does not vanish$.$

\subsection{Pigeonhole principle\label{PH}}

We can split $\mathcal{D}_{n}$ into $\mathcal{D}_{n}(p)$'s where $\mathcal{D}%
_{n}(p)$ is the set of $c$'s such that $p$ divides $a(c).$ But $(\ref{bound}%
) $ implies that no more than $n^{r}$ primes bigger than $e^{Cn}$
can divide a single $a(c).$ In particular for every finite set
$\mathcal{P}$ of primes larger that $e^{Cn}$ we have
\begin{equation*}
|\mathcal{P}|\cdot \min_{p\in \mathcal{P}}|\mathcal{D}_{n}(p)|\leq
\sum_{p\in \mathcal{P}}|\mathcal{D}_{n}(p)|\leq n^{r}\cdot |\cup \mathcal{D}%
_{n}(p)|\leq n^{r}|\mathcal{D}_{n}|
\end{equation*}%
Thus, given $\varepsilon >0,$ if $|\mathcal{P}|\geq |\mathcal{D}%
_{n}|^{2\varepsilon },$ there must be a prime $p\in \mathcal{P}$ with $|%
\mathcal{D}_{n}(p)|<|\mathcal{D}_{n}|^{1-\varepsilon }.$

Now we pass to the second part of the proof:\ namely it remains to
show that
if a prime $p$ satisfies $|\mathcal{D}_{n}(p)|\leq |\mathcal{D}%
_{n}|^{1-\varepsilon }$ then a lower bound on $g(p)$ can be
deduced. This is of course the place where we will use the results
of \cite{BG}.

\subsection{Modified \protect\cite{BG}: subgroup non-concentration implies
gap}

In fact, rather than the main statement of \cite{BG}, which gave a
lower bound on $g(a,b)$ in terms of the girth of the pair $(a,b)$,
we are going to explain how the proof of \cite{BG} allows to
obtain a similar lower bound
out of the weaker hypothesis that the simple random walk on $<a,b>$ at time $%
constant\times log(p)$ gives a weight of at most $1/p^{constant}$
to every
proper subgroup $\mathrm{SL}(2,p).$ Namely, writing $\mu _{(a,b)}=\frac{1}{4}%
(\delta _{a}+\delta _{a^{-1}}+\delta _{b}+\delta _{b^{-1}}),$

\begin{theorem}[Modified \protect\cite{BG}]
\label{BGthm}There is a function $\delta =\delta (\tau ,\gamma
)>0$ such that for every $\tau ,\gamma >0$, every large prime $p,$
and every generating pair $(a,b)$ in $\mathrm{SL}(2,p)$ such that
\begin{equation}
\sup_{H}\mu _{(a,b)}^{(\tau \log _{3}p)}(H)\leq p^{-\gamma }
\label{wg}
\end{equation}%
(where the sup is taken over all proper subgroups $H$ of
$SL(2,p)$), we have $g(a,b)>\delta .$
\end{theorem}

Let $2l=\tau \log _{3}p$ and let $\nu =\mu _{(a,b)}^{(l)}$. In
\cite{BG} the logarithmic girth condition is used to verify that
$\nu $ satisfies the two conditions of the $l^{2}$ flattening
lemma (Proposition 2 in \cite{BG}), namely: (a) $\Vert \nu \Vert
_{\infty }<p^{-\gamma }$ and (b) $\nu ^{(2)}(H)<p^{-\gamma }$ for
all proper subgroups $H$ of $SL(2,p)$. Condition (b) follows
immediately from \eqref{wg}; condition (a) also easily follows
by applying \eqref{wg} with the trivial subgroup $H=\{e\}$. Indeed, we have $%
\mu _{(a,b)}^{(2l)}(e)=\Vert \mu _{(a,b)}^{(l)}\Vert _{2}^{2}$ and
$\Vert \mu _{(a,b)}^{(l)}\Vert _{\infty }\leq \Vert \mu
_{(a,b)}^{(l)}\Vert _{2}$, thus we obtain $\Vert \mu
_{(a,b)}^{(l)}\Vert _{\infty }<p^{-\gamma }$.

\subsection{Proof of Theorem \protect\ref{mainbis}}

Since $|\mathcal{C}_{n}|=4^{4Kn}$ is an exponential function of
$n,$ observe
that it is enough to prove the theorem for $X$ of the form $X=|\mathcal{C}%
_{n}|$ for some $n.$ If there are less that $X^{2\varepsilon }$
primes between $X$ and $X^{A}$ there is nothing to prove.
Otherwise \S \ref{PH}
implies that there is a prime $p$ between $X$ and $X^{A}$ such that $|%
\mathcal{D}_{n}(p)|<|\mathcal{D}_{n}|^{1-\varepsilon }$. Let
$(a,b)\in \mathrm{SL}(2,p)^{2}$ a generating pair. Note that Lemma
\ref{gen} implies
that for every $c\in \mathcal{D}_{n}\backslash \mathcal{D}_{n}(p)$ there is $%
w_{1},w_{2}$ in $B(N)^{2}$ such that the corresponding commutator of length $%
4n$ in $w_{1},w_{2}$ is not $1$. Suppose that for every $w_{1},w_{2}$ in $%
B(N)^{2}$ there is a proper subgroup $H$ of $\mathrm{SL}(2,p)$ such that $%
\mu _{(w_{1},w_{2})}^{\ast n}(H)\geq 4^{-\varepsilon n/2}.$ Recall
the subgroup structure of $\mathrm{SL}(2,p)$ : every proper
subgroup either has
cardinality at most $60$, or must be solvable of solvability class at most $%
2 $ (see \cite{Dick}). Let $S_{n}^{(i)}$ be four independent
simple random
walks starting at $1$ on the subgroup of $\mathrm{SL}(2,p)$ generated by $%
(w_{1},w_{2}).$ Then if $H$ has cardinality at most $60,$ $\mu
_{(w_{1},w_{2})}^{\ast n}(H)\leq 60\cdot
\mathbb{P}(S_{n}^{(1)}=1),$ while
if $H$ is solvable, $\mu _{(w_{1},w_{2})}^{\ast n}(H)^{4}\leq \mathbb{P}%
([[S_{n}^{(1)},S_{n}^{(2)}],[S_{n}^{(3)},S_{n}^{(4)}]]=1)$. It
follows in both cases that are at least $4^{4(1-\varepsilon /2)n}$
quadruples of paths of length $n$ whose $4$-fold commutator
vanishes. Then we can count at least $4^{4(1-\varepsilon
/2)nK}=|\mathcal{C}_{n}|^{1-\varepsilon /2}$ assignments $c\in
\mathcal{C}_{n}$ for which all corresponding words vanish at
$\left( a,b\right) $ ; hence at least
$|\mathcal{C}_{n}|^{1-\varepsilon }$ assignments $c\in $
$\mathcal{D}_{n}$ for which all corresponding commutator words
vanish at $\left( a,b\right) $ (recall $|\mathcal{C}_{n}\backslash
\mathcal{D}_{n}|\leq |\mathcal{C}_{n}|^{1-\alpha }$ for some explicit $%
\alpha >\varepsilon /2>0$ by Corollary \ref{c vs d}). By the
preceding
remark, those must belong to $\mathcal{D}_{n}(p).$ Hence $|\mathcal{D}%
_{n}(p)|\geq |\mathcal{C}_{n}|^{1-\varepsilon }$ a contradiction.
Therefore, there must exist $w_{1},w_{2}$ such that $\mu
_{(w_{1},w_{2})}^{\ast n}(H)\leq 4^{-\varepsilon n/2}$ for all
proper subgroups $H$ in $SL(2,p)$.
But $p\in \lbrack X,X^{A}]$ and $X=4^{Kn}.$ Hence we may apply Theorem \ref%
{BGthm} and deduce that $g(w_{1},w_{2})>\delta .$ But this readily
implies that $g(a,b)>\delta /N$, which ends the proof of Theorem
\ref{mainbis}.

\subsection{Proof of Corollary \protect\ref{c0}}

For a Cayley graph of a finite group $G$ generated by a symmetric
set $S$ containing $1$ to be a $c$-expander, it is sufficient that
for every $f\in \ell _{0}^{2}(G)$ (functions with zero average on
$G$) there is an $s\in S$
such that $||s\cdot f-f||_{2}\geq 2\sqrt{c}\cdot ||f||_{2}.$ Indeed take $f=a%
\mathbf{1}_{A}-b\mathbf{1}_{A^{c}}$ with a choice of $a$ and $b$ such that $%
a|A|=b|A^{c}|$. Then, since $|A^{c}|\geq |A|,$ we have $a\geq b$ and $%
||f||_{2}^{2}=(a+b)a|A|\geq \frac{(a+b)^{2}}{2}|A|,$ while
$||s\cdot f-f||_{2}^{2}=(a+b)^{2}|sA\Delta A|.$ Thus $|\partial
A|\geq |sA\backslash A|=\frac{1}{2}|sA\Delta A|\geq c\cdot |A|$.

Now observe (by the triangle inequality) that if there is a
constant $N$ such that for any $f\in \ell _{0}^{2}(G)$ there is
$\gamma \in S^{N}$ such that $||\gamma \cdot f-f||_{2}\geq
2N\sqrt{c}\cdot ||f||_{2},$ then there must also be some $s\in S$
such that $||s\cdot f-f||_{2}\geq 2\sqrt{c}\cdot ||f||_{2}$.
Therefore Corollary \ref{c0} will follow from Theorem \ref{main}
if we can show that there is a constant $N$ independent of $S$ and
$p$ such that $S^{N}$ contains two elements $\{a,b\}$ that
generate $\mathrm{SL}(2,p). $ As in the last paragraph, observe
that every proper subgroup either has cardinality at most $60$, or
must be solvable of solvability class at most $2
$, therefore for $\{a,b\}$ to generate $\mathrm{SL}(2,p)$ it is enough that $%
a$ and $b$ have no relation of length at most 60 say. The
existence of such a constant $N$ is an immediate consequence of
the strong Tits alternative, i.e. Theorem \ref{ST} (see \cite{B2}
for more details on this derivation). This ends the proof of
Corollary \ref{c0}.

\textbf{Acknowledgement:} We are grateful to H. Helfgott for
pointing out \cite{Ho}.

\end{document}